\documentclass[11pt]{article}
\usepackage{amssymb,amsmath}
\usepackage[dvips]{graphicx}
\usepackage[latin1]{inputenc}
\newtheorem{dfn}{Definition}[section]

\newtheorem{thm}{Theorem}[section]
\newtheorem{thm*}{Theorem}
\newtheorem{cor}{Corollary}[section]
\newtheorem{prop}{Proposition}[section]

\newcommand{\Pf}{{\em Proof}. }
\newcommand{\EPf}{\begin{flushright} $\Box$ \end{flushright}}

\newcommand{\dif}{{\rm d}}
\newcommand{\ca}{{\rm ca}}

\newcommand{\e}{\mathbf{e}}

\begin{document}
\title{Gauss-Bonnet  theorem in sub-Riemannian Heisenberg space $\mathbb H^1$}
\author{ Jos\'e M. M. Veloso
         \\ Instituto de Ci\^encias Exatas e Naturais
         \\ Universidade Federal do Par\'a
                  \\ veloso@ufpa.br
       \and Marcos M. Diniz
         \\ Instituto de Ci\^encias Exatas e Naturais
         \\ Universidade Federal do Par\'a
         \\ mdiniz@ufpa.br}
\maketitle
\begin{abstract}
We  prove a version of Gauss-Bonnet theorem in sub-Riemannian Heisenberg space $\mathbb H^1$.  The sub-Riemannian distance makes  $\mathbb H^1$ in a metric space and consenquently with a spherical Hausdorff measure. Using this measure, we define a Gaussian curvature at points of a surface $S$ where the sub-Riemannian distribution is transverse to the tangent space of $S$. If all points of $S$ have this property, we prove a Gauss-Bonnet formula and for compact surfaces (which are topologically a torus) we obtain $\int_SK=0$. 
\end{abstract}

{\bf Primary subject:} 53C17

\section{Introduction} In this paper, we prove a Gauss-Bonnet type theorem for surfaces inside Heisenberg group $\mathbb H^1$. In this space consider a distribution $D$ generated by vector fields $$\e_{1}=\frac{\partial}{\partial x}-\frac 1 2y\frac{\partial}{\partial z} \ \ \ ;\ \ \  \e_2=\frac{\partial}{\partial y}+\frac 1 2x\frac{\partial}{\partial z},$$ and a scalar product in $D$ such that $\e_1,\e_2$ are orthonormal. Complete these vector fields to a basis of left invariant vector fields in $\mathbb H^1$  introducing 
$$\e_0=[\e_1,\e_2]=\frac{\partial}{\partial z}.$$ 
Therefore, if $\e^0, \e^1, \e^2$ are dual forms to $\e_0,\e_1,\e_2$, then the volume element invariant by the group action is $dV=\e^0\wedge \e^1\wedge \e^2$.

With the scalar product in $D$,  consider the distance between two points as the infimum of length of curves tangent to $D$ that connect them. With this distance, $\mathbb H^1$ is a metric space with Hausdorff dimension four and the differentiable surfaces have dimension three.
At points of a surface $S$ where the distribution $D$ does not coincide with $TS$, the intersection $D\cap TS$ has dimension one, and we obtain a direction  called \emph{characteristic} at this point of $S$. \emph{We suppose every point of surface $S$ has a characteristic direction.} The vector field  \emph{normal horizontal} $\eta$ is an unitary vector field in $D$ and orthogonal to the characteristic direction which we suppose is globally defined. Given a compact set  $K\subset S$, the 3-dimensional (spherical) Hausdorff measure of $K$ is given by $\int_Ki(\eta)dV$. A curve  transverse to $D$ has Hausdorff dimension two and its (spherical) Hausdorff measure is given by $\int_{\gamma}\e^0$. For more details, see  \cite{Gromov.1996,Montgomery.2002,Pansu.1982,Pansu.1982a}.

To prove a Gauss-Bonnet theorem,  we need a concept of curvature of surfaces. The image by the left transport of the normal horizontal  in a neighborhood of a point in $S$ is contained in $S^1\subset T_0\mathbb H^1$, therefore the normal horizontal does not suit as a Gauss map. But we can consider the 1-form   $\eta^{*}$ defined on $S$ by $\eta^ *(\eta)=1$ and $\eta^*|_{TS}=0$. The  analogous of Gauss application is 
$$\begin{array}{rccc}
g:=\exp\circ L^*\circ\eta^*:&S&\rightarrow&\mathbb H^1\\
&p&\mapsto&\exp(L_p^*(\eta^*(p))),
\end{array}$$
with image  in the cilinder $S^1\times \mathbb R$. Then we define  
\begin{equation}\label{curvg}
K(p)=\lim_{U\rightarrow\{p\}}\frac{\int_{g(U)}i(\tilde\eta)dV}{\int_{U}i(\eta)dV}
\end{equation} 
as the \emph{Gaussian curvature} of surface $S$ at point $p$, where $\tilde\eta$ is the horizontal normal to $g(S)$. 

Consider the \emph{adapted  covariant derivative} $\overline \nabla$, defined in \cite{Falbel.1993}, for which the left invariant vector fields are parallel.  We then define a covariant derivative $\nabla$ on $S$ by projecting $\overline\nabla$ in the direction of $\eta^*$,
$$
\nabla_XY=\overline\nabla_XY-\eta^*(\overline\nabla_XY)\eta,
$$
where $X,Y$ are vector fields on $S$. The relevant fact is that the curvature associated to $\nabla$ coincides with the one defined by Gauss map  (\ref{curvg}).

To get the local form of Gauss-Bonnet theorem, we still need the concept of geodesic curvature for curves in the surface. We consider curves transverse to characteristic directions, and for these curves we define the tangent field  $$T=\frac{\gamma'}{\e^0(\gamma')}.$$ 
If $N$ is an unitary field in the characteristic direction along the transverse curve $\gamma$, with orientation conveniently chosen, then  we have $\nabla_TT=kN$, and  $k$ is the \emph{curvature} of $\gamma$. Finally, to characterize the  variation of directions of two transverse curves by a same vertex, we define the \emph{corner area} between two tangent vectors of $S$ in a point by $$\mbox{ca}(v,w)=\frac{dV(\eta,v,w)}{\e^0(v)\e^0(w)}.$$

With these preliminaries, we state

\begin{thm}(Gauss-Bonnet formula) Let ${R}$ be a region contained in a coordinate domain $U$ of $S$ such that $T_pS \neq D_p$, for all $p\in U$, let the bounding curve $\gamma$ of $R$ be a simple closed transverse curve, and let $ca_1,\ldots,ca_r$ be the exterior corner areas of $\gamma$. Then
$$\int_{\gamma}k+\sum_{j=1}^rca_j+\int_RK=0,$$
where $k$ is the curvature function on $\gamma$ and $K$ is the  Gaussian curvature function on $R$.
\end{thm}

If the surface $S$ is compact and oriented, then there exists a characteristic no-null vector field on $S$, therefore $S$ is diffeomorphic to a torus. In this case, we obtain the corollary:

\begin{cor}
Suppose $S$ is a differentiable compact surface in $\mathbb H^1$ such that $T_pS \neq D_p$, for all $p\in S$. Then
$$\int_SK=0.$$
\end{cor}

\section{The Heisenberg group}
We denote by $\mathbb H^1$ the Heisenberg nilpotent Lie group whose manifold is $\mathbb R^3$, with Lie algebra ${H}^1=V_1\oplus V_2$, $\dim V_1=2$, $\dim V_2=1$, and  
$$[V_1,V_1]=V_2\ \ \ ;\ \ \  [V_1,V_2]=[V_2,V_2]=0.$$
Since   $\mathbb H^1$ is nilpotent, the exponential map $\exp:H^1\rightarrow \mathbb H^1$ is a diffeomorphism. Let be $\e_1,\e_2$  a basis of $V_1$ and $\e_0=[\e_1,\e_2]\in V_2$. By applying the Baker-Campbell-Hausdorff formula we have
$$\exp^{-1}(\exp(X)\exp(Y))=X+Y+\frac{1}{2}[X,Y].$$
Since $[\e_1,\e_2]=\e_0$, writing $X=x_1\e_1+y_1\e_2+z_1\e_0$, and $Y=x_2\e_1+y_2\e_2+z_2\e_0$, we get 
$$ X+Y+\frac{1}{2}[X,Y]=(x_1+x_2)\e_1+(y_1+y_2)\e_2+(z_1+z_2+\frac{1}{2}(x_1y_2-x_2y_1))\e_0. $$ 
We identify $\mathbb H^1$ with $\mathbb R^3$ by identifying $(x,y,z)$ with $\exp(x\e_1+y\e_2+z\e_0)$, and this  is known as  canonical coordinates of first kind or exponential coordinates. In these coordinates, the group operation is
$$(x_1,y_1,z_1)(x_2,y_2,z_2)=(x_1+x_2,y_1+y_2,z_1+z_2+\frac{1}{2}(x_1y_2-x_2y_1)),$$
the exponential is 
$$\exp(x\e_1+y\e_2+z\e_0)=(x,y,z),$$
and the left invariant vector fields $\e_1,\e_2,\e_0$ are given by
\[\left\{
\begin{aligned}
\e_1&=\frac{\partial}{\partial x}-\frac 1 2y\frac{\partial}{\partial z},\\[.2cm]
\e_2&=\frac{\partial}{\partial y}+\frac 1 2x\frac{\partial}{\partial z},\\[.2cm]
\e_0&=\frac{\partial}{\partial z},
\end{aligned}
\right.
\]
with brackets 
$[\e_1,\e_2]=\e_0,$ and $[\e_0,\e_1]=[\e_0,\e_2]=0$.
 The dual basis is
\[
\left\{
\begin{aligned}
\e^1&=\dif x,\\
\e^2&=\dif y,\\
\e^0&=\dif z+\frac 1 2(y\dif x-x\dif y),
\end{aligned}
\right.
\]
with $\dif \e^0=-\e^1\wedge \e^2$, $\dif \e^1=\dif \e^2=0$. For more details, see \cite{Capogna.2007}.

We identify naturally $T\mathbb H^1$ with $T^*\mathbb H^1$ by
$a\e_1+b\e_2+c\e_0$ with $a\e^1+b\e^2+c\e^0$, and through this identification we identify $H^1$
with $(H^1)^*$. Therefore we can define the exponential map on the dual by
$$\begin{array}{rccc}
\exp:&(H^1)^*&\rightarrow&\mathbb H^1\\
&x\e^1+y\e^2+z\e^0&\mapsto&(x,y,z) .
\end{array}$$
The left translation is defined by
$$L_{(x,y,z)}(x_1,y_1,z_1)=(x,y,z)(x_1,y_1,z_1),$$
and
$$L_{(x,y,z)}^{-1}=L_{(-x,-y,-z)}.$$
Let be $D\subset T\mathbb H^1$ the two-dimensional distribution generated by the vector fields $\e_1,\e_2$, so that $D$ is the null space of $\e^0$. On $D$, we define a scalar product  $ \langle\, , \rangle$, such that $\{\e_1,\e_2\}$ is an orthonormal basis of $D$. An operator $J:D\rightarrow D$ is well-defined by
$$J(a\e_1+b\e_2)=-b\e_1+a\e_2.$$

The \emph{element of volume} $\dif V$ in $\mathbb H^1$ is $\dif V=\e^0\wedge \e^1\wedge \e^2=\dif x\wedge \dif y\wedge \dif z.$ A differentiable curve  $\gamma:[a,b]\subset \mathbb R\rightarrow \mathbb H^1$ is \emph{transversal} if $\e^0(\gamma'(t))\neq 0$, for every $t\in [a,b]$. We say that a transversal curve $\gamma$ is \emph{unitarily parametrized} if $|\e^0(\gamma'(t))|=1$, for every $t\in [a,b]$.

\section{The adapted covariant derivative}

If $X,Y$ are vector fields on $\mathbb H^1$, we define the adapted covariant derivative introduced in \cite{Falbel.1993} by:
$$\overline\nabla_XY=\sum_{i=0}^2 \dif b_j(X)\e_j,$$
where $Y=b_0\e_0+b_1\e_1+b_2\e_2$. Then $\overline\nabla$ is null at left invariant vector fields on $\mathbb H^1$.

\begin{prop}The covariant derivative 
$\overline\nabla$ has the following properties:
\begin{enumerate}
  \item
If $Y\in\underline{D}$, then $\overline\nabla_XY\in \underline{D}$ for all $X\in \underline{T\mathbb H^1}$;
\item
 If $Y,Z\in\underline{D}$, then 
$$\overline\nabla_X\langle Y,Z\rangle =\langle\overline\nabla_XY,Z\rangle +\langle Y,\overline\nabla_XZ\rangle ,$$
 for all $X\in \underline{T\mathbb H^1}$;
 \item The torsion $\overline{T}$ of $\overline\nabla$ is
 $$\overline{T}=-\e^1\wedge \e^2\otimes \e_0=\dif \e^0\otimes \e_0,$$
 \item The curvature $\overline{K}$ of $\overline\nabla$ is null.
\end{enumerate}
\end{prop}

\Pf We shall proceed the proof of \emph{3}, the others being similar. If $X=\sum_{i=0}^2a_j\e_j$ and $Y=\sum_{i=0}^2b_j\e_j$, then
$$\overline{T}(X,Y)=\overline\nabla_XY-\overline\nabla_YX-[X,Y]=-(a_1b_2-a_2b_1)e_0
=-\e^1\wedge \e^2(X,Y)\e_0=\dif\e^0(X,Y)\e_0.
$$
\EPf

Observe that the covariant derivative in the cotangent bundle $(T\mathbb H^1)^*$ satisfies $\overline\nabla \e^i=0$, for $i=0,1,2$.

\section{Surfaces in $\mathbb H^1$}
Suppose $S$ is an oriented differentiable two-dimensional manifold in $\mathbb H^1$. Note that $\dim(D\cap TS)\geq 1$, and, since $\dif \e^0=\e^1\wedge \e^2$, the set of points where the tangent space of $S$ coincides with the distribution has empty interior. We denote by $\Sigma$ this set and by $S'$ its complement on $S$,
\[\Sigma=\{x\in S:  \dim(D_x\cap T_xS)=2\}\ \ \ ;\ \ \ S'=S-\Sigma.\] 
The set $S'$ is open in $S$. In what follows we will suppose $\Sigma=\emptyset$, so  $S=S'$. With this hypothesis on $S$, the one-dimensional vector subbundle $D\cap TS$ is well defined. Suppose $U\subset S$ is an open set such that we can define a unitary vector field $f_1$ with values in $D\cap TS$, so $\langle f_1,f_1\rangle =1$. 
\begin{dfn}The unitary vector field $\eta\in \underline{D}$ defined by
$$\eta=-Jf_1$$
is the \emph{horizontal normal} to $S$.
\end{dfn}
Then we can define $\eta^*\in (T\mathbb H^1)^*|_S$ by 
\[\eta^*(\eta)=1 \ \ \ ;\ \ \ \eta(TS)=0.\]
 We call $\eta^*$ the \emph{horizontal conormal} to $S$. 

\begin{dfn}The application
$$\begin{array}{rccc}
g:=\exp\circ L^*\circ\eta^* :&S&\rightarrow&\mathbb H^1\\
&p&\mapsto &\exp(L_p^*(\eta^*(p)))
\end{array}
$$
is the \emph{Gauss map} of $S$.
\end{dfn}
Let be 
$$f_2=\e_0-\eta^*(\e_0)\eta.$$
Then $\{f_1,f_2\}$ is a \emph{special} basis of $TS$ on the open set $U$. If
$$\eta=\cos\alpha \e_1+\sin\alpha \e_2,$$
for some real function $\alpha$ on $U$, reducing $U$ if necessary, then
$$f_1=-\sin\alpha \e_1+\cos\alpha \e_2,$$
and, if we denote by $A=-\eta^*(\e_0)$, we write
$$f_2=\e_0+A\eta.$$
The dual basis of $(T\mathbb H^1)^*$ on $S$ is 
\[
\left\{
\begin{aligned}
\eta^*&=\cos\alpha \e^1+\sin\alpha \e^2-A\e^0,\\
f^1&=-\sin\alpha \e^1+\cos\alpha \e^2,\\
f^2&=\e^0.
\end{aligned}
\right.
\]
The inverse relations are
\[
\left\{
\begin{aligned}
\e^0&=f^2,\\
\e^1&=\cos\alpha\,\eta^*-\sin\alpha f^1+A\cos\alpha f^2,\\
\e^2&=\sin\alpha\,\eta^*+\cos\alpha f^1+A\sin\alpha f^2,
\end{aligned}
\right.
\]
and 
$$\e^1\wedge \e^2=\eta^*\wedge f^1-Af^1\wedge f^2.$$
Also, it follows
\[
\left\{
\begin{aligned}
\dif f^1&=-\dif\alpha\wedge\eta^*-A\,\dif\alpha\wedge f^2,\\
\dif f^2&=-\eta^*\wedge f^1+Af^1\wedge f^2,\\
\dif\eta^*&=(\dif\alpha+A^2f^2+A\eta^*)\wedge f^1-\dif A\wedge f^2,
\end{aligned}
\right.
\]
and, since $\eta^*=0$ on $S$, we get
\[
\left\{
\begin{aligned}
\dif f^1&=-A\,\dif\alpha\wedge f^2,\\
\dif f^2&=Af^1\wedge f^2,\\
0&=(\dif\alpha+A^2f^2)\wedge f^1-\dif A\wedge f^2.
\end{aligned}
\right.
\]
From this last relation, we obtain 
$$\dif \alpha(f_2)=-(\dif A(f_1)+A^2).$$

\begin{dfn}The \emph{element of area} in $S$ is
$$i(\eta)\dif V.$$
\end{dfn}

Since $\dif V=\eta^*\wedge f^1\wedge f^2$, then $\dif S=f^1\wedge f^2$.
Let's find the area  of $g(R)$ for a region $R\subset S$. Observe that, for all $p\in S$, 
$$g(p)=(\cos\alpha(p),\sin\alpha(p) ,-A(p)).$$
Then $g(R)$ is contained on the cylinder $C=\{(x,y,z):x^2+y^2=1\}.$ The tangent space $TC$ is generated by 
\[
\left\{
\begin{array}{ll}
-y\e_1+x\e_2+\dfrac 1 2 \e_0 & (=-y\dfrac{\partial}{\partial x}+x\dfrac{\partial}{\partial y})\\
\ \e_0 &(=\dfrac{\partial}{\partial z}).
\end{array}
\right.
\]
It follows that 
\[
\left\{
\begin{aligned}
\widetilde{f}_1&=-y\e_1+x\e_2\\
\widetilde{f}_2&=\e_0,
\end{aligned}
\right.
\]
so
$$\widetilde{\eta}=-J(\widetilde{f}_1)=x\e_1+y\e_2.$$
The element of area on $C$ is $\dif\widetilde{S}=\widetilde{f}^1\wedge \widetilde{f}^2$. Then
$$
\begin{aligned}
\mathrm{Area}(g(R))&=\int_{g(R)}\dif\widetilde{S}=\int_{g(R)}\widetilde{f}^1\wedge \widetilde{f}^2=\int_Rg^*(\widetilde{f}^1\wedge \widetilde{f}^2)=\int_Rg^*(\widetilde{f}^1\wedge \widetilde{f}^2)(f_1,f_2)f^1\wedge f^2\vspace{.3cm}\\
&=\int_R(\widetilde{f}^1\wedge \widetilde{f}^2)(g_*f_1,g_*f_2)\dif{S}.\\
\end{aligned}
$$
Now, 
$$
\begin{aligned}
\dif g&=-\sin\alpha \dif \alpha\otimes \frac{\partial}{\partial x}+\cos\alpha \dif \alpha\otimes \dfrac{\partial}{\partial y}-\dif A\otimes \dfrac{\partial}{\partial z}\vspace{.3cm}\\
&=-\sin\alpha \dif \alpha\otimes (\e_1+\dfrac 1 2 \sin \alpha \e_0)+\cos\alpha \dif \alpha\otimes (\e_2-\dfrac 1 2 \cos \alpha \e_0)-\dif A\otimes \e_0\vspace{.3cm}\\
&=\dif \alpha\otimes \widetilde{f}_1-(\dfrac 1 2 \dif \alpha+\dif A)\otimes \widetilde{f}_2,
\end{aligned}
$$
and so
$$
\begin{aligned}
(\widetilde{f}^1\wedge \widetilde{f}^2)(g_*f_1,g_*f_2)&=-\dif \alpha(f_1)\left(\frac 1 2 \dif \alpha(f_2)+\dif A(f_2)\right)+\dif \alpha(f_2)\left(\frac 1 2 \dif \alpha(f_1)+\dif A(f_1)\right)\\
&=-\dif \alpha\wedge \dif A(f_1,f_2).
\end{aligned}$$
We just proved that
$$\mathrm{Area}(g(R))=\int_R -\dif \alpha\wedge \dif A(f_1,f_2) \dif S.$$
As $\mathrm{Area}(R)=\int_R  \dif S$, we obtain from (\ref{curvg}) that
$$K=-\dif \alpha\wedge \dif A(f_1,f_2).$$
\begin{prop}The Gaussian curvature $K$ of $S$ is given by
$$K=-\dif \alpha\wedge \dif A(f_1,f_2).$$
\end{prop}

\section{The projection of $\overline{\nabla}$ by $\eta^*$}

Given  $X,Y\in \underline{TS}$, we define 
$$\nabla_XY=\overline{\nabla}_XY-\eta^*(\overline{\nabla}_XY)\eta.$$

\begin{prop}\label{pnabla}The operator $\nabla$ is a covariant derivative in $TS$, and satisfies:
\begin{enumerate}
\item $\nabla f_1=0$;
\item $\nabla f_2=A\,\dif \alpha\otimes f_1$;
\item $\nabla f^1=-A\,\dif \alpha\otimes f^2$;
\item $\nabla f^2=0.$
\end{enumerate}
\end{prop}
\Pf It is clear that, if $X,Y\in TS$, then $\nabla_XY\in TS$, $\nabla_XY$ is linear on $X$ and additive on $Y$. Furthermore, if $f$ is a real function on $S$, we have 
$$\nabla _XfY=
\dif f(X)Y+f\overline{\nabla}_XY-\eta^*(\dif f(X)Y+f\overline{\nabla}_XY)\eta=\dif f(X)Y+f\nabla_XY.$$
Finally, 
\begin{itemize}
\item[\emph{1}.] $\overline{\nabla}_Xf_1=\overline{\nabla}_X(-\sin\alpha \e_1+\cos\alpha \e_2)=\dif \alpha(X)(-\cos\alpha \e_1-\sin\alpha \e_2)=-\dif \alpha(X)\eta$, so $\nabla_Xf_1=0$.
\item[\emph{2}.] $\overline{\nabla}_Xf_2=\overline{\nabla}_X(\e_0+A\eta)=\dif A(X)\eta+A\overline{\nabla}_X(\cos\alpha \e_1+\sin\alpha \e_2)=\dif A(X)\eta+A\,\dif \alpha(X)(-\sin\alpha \e_1+\cos\alpha \e_2)=\dif A(X)\eta+A\,\dif \alpha(X)f_1$, so $\nabla_Xf_2=A\,\dif \alpha(X)f_1$.
\item[\emph{3}.] $(\nabla_Xf^1)(f_1)=-f^1(\nabla_Xf_1)=0$ and $(\nabla_Xf^1)(f_2)=-f^1(\nabla_Xf_2)=-A\,\dif \alpha(X)$ so $\nabla_Xf^1=-A\,\dif \alpha(X)f^2$.
\item[\emph{4}.] $(\nabla_Xf^2)(f_1)=-f^2(\nabla_Xf_1)=0$ and $(\nabla_Xf^2)(f_2)=-f^2(\nabla_Xf_2)=0$ so $\nabla_Xf^2=0$.
\end{itemize}
\EPf

It follows from this proof that $\overline{\nabla}_X\eta=\nabla_X\eta=\dif \alpha(X)f_1$ and $\overline{\nabla}_X\eta^*=\dif \alpha(X)f^1-\dif A(X)f^2$, for $X\in TS$.

\begin{dfn}
The covariant derivative $\nabla$ is the \emph{adapted} covariant derivative on $S$.
\end{dfn}

\begin{prop}
The torsion $T$ of $\nabla$ is $T=Af^1\wedge f^2\otimes f_2$.
\end{prop}
\Pf We have 
$$T(X,Y)=\overline\nabla_XY-\overline\nabla_YX-[X,Y]-\eta^*(\overline\nabla_XY-\overline\nabla_YX-[X,Y])\eta=\overline{T}(X,Y)-\eta^*(\overline{T}(X,Y))\eta,$$
so 
$$T=-\e^1\wedge \e^2\otimes ( \e_0-\eta^*( \e_0)\eta)=Af^1\wedge f^2\otimes f_2.$$

\begin{prop}
The curvature tensor $R$ of $\nabla$ is $R=\dif A\wedge d\alpha\otimes f^2\otimes f_1$.
\end{prop}
\Pf Clearly $R(X,Y)f_1=0$, and 
$$
\begin{aligned}
R(X,Y)f_2&=\nabla_X\nabla_Yf_2-\nabla_Y\nabla_Xf_2-\nabla_{[X,Y]}f_2 \vspace{.3cm}\\ 
&=\nabla_X(A \dif \alpha(Y)f_1)-\nabla_Y(A\dif \alpha(X)f_1)-A\dif \alpha([X,Y])f_1\vspace{.3cm}\\
&=\left(X(A \dif \alpha(Y))-Y(A\dif \alpha(X))-A\dif \alpha([X,Y])\right)f_1\vspace{.3cm}\\
&=\dif (A\dif \alpha)(X,Y)f_1.\vspace{.3cm}\\
\end{aligned}
$$
\EPf
\begin{prop}
The Gaussian curvature $K$ is given by
\begin{equation}\label{K}
K=\langle R(f_1,f_2)f_2,f_1\rangle= \dif A\wedge \dif \alpha(f_1,f_2).
\end{equation}
\end{prop}

\section{The second fundamental form}
From the equation
$$\overline{\nabla}_XY=\nabla_XY+\eta^*(\overline{\nabla}_XY)\eta=\nabla_XY-(\overline{\nabla}_X\eta^*)(Y)\eta,$$
for $X,Y\in \underline{TS}$, we define a bilinear form $V\!:TS\times TS\rightarrow \mathbb R$:
\begin{dfn}
The bilinear form  $V\!:TS\times TS\rightarrow \mathbb R$, defined by
$$V(X, Y)=-(\overline{\nabla}_X\eta^*)(Y)$$
is the \emph{second fundamental form} associated to $S$.
\end{dfn}
From
$$\overline{\nabla}\eta^*=\overline{\nabla}(\cos\alpha \e^1+\sin\alpha \e^2-A\e^0)=\dif \alpha\otimes(-\sin\alpha \e^1+\cos\alpha \e^2)-\dif A\otimes \e^0,$$
we get
$$V(X,Y)=-\dif \alpha(X) f^1(Y)+\dif A(X) f^2(Y).$$
The second fundamental form is not symmetric in general. In fact, for $X,Y\in TS$,
$$\begin{aligned}
V(X,Y)-V(Y,X)&=-(\overline{\nabla}_X\eta^*)(Y)+(\overline{\nabla}_Y\eta^*)(X)=\eta^*(\overline{\nabla}_XY)-\eta^*(\overline{\nabla}_YX)\\
&=\eta^*(\overline{T}(X,Y))=\dif \e^0(X,Y)\eta^*(\e_0)=-A \dif f^2(X,Y)\\
&=-A^2f^1\wedge f^2(X,Y).
\end{aligned}$$

\begin{thm}The curvature $K$ and the second fundamental form $V$ satisfy:
\begin{enumerate}
\item  (Gauss equation) $K(X,Y)Z=(-\dif \alpha(X)V(Y,Z)+\dif \alpha(Y)V(X,Z))f_1;$
\item (Codazzi equation) $\nabla_XV(Y,Z)-\nabla_YV(X,Z)+V(T(X,Y),Z)=0. $
\end{enumerate}
\end{thm}
\Pf By applying the definition of curvature, we obtain
$$\begin{aligned}
K(X,Y)Z=&\overline{\nabla}_X(\overline{\nabla}_YZ-V(Y,Z)\eta)-V(X,\nabla_YZ)\eta-\overline{\nabla}_Y(\overline{\nabla}_XZ-V(X,Z)\eta)\\
&+V(Y,\nabla_XZ)\eta-\overline{\nabla}_{[X,Y]}Z+V([X,Y],Z)\eta\\[.2cm]
=&\overline{K}(X,Y)Z-X(V(Y,Z))\eta-V(Y,Z)\dif \alpha(X)f_1+Y(V(X,Z))\eta\\
&+V(X,Z)\dif \alpha(Y)f_1-V(X,\nabla_YZ)\eta+V(Y,\nabla_XZ)\eta+V([X,Y],Z)\eta\\[.2cm]
=&\left(-\nabla_XV(Y,Z)-V(\nabla_XY,Z) +\nabla_YV(X,Z)+V(\nabla_YX,Z)+V([X,Y],Z)\right)\eta\\
&+ \left(-V(Y,Z)\dif \alpha(X)+V(X,Z)\dif \alpha(Y)\right)f_1\\[.2cm]
=&-\left(\nabla_XV(Y,Z)-\nabla_YV(X,Z)+V(\nabla_XY-\nabla_YX-[X,Y],Z)\right)\eta\\
&-\left(V(Y,Z)\dif \alpha(X)-V(X,Z)\dif \alpha(Y)\right)f_1.
\end{aligned}
$$
Since $K(X,Y)Z\in TS$, we obtain
$$K(X,Y)Z=-\left(V(Y,Z)\dif \alpha(X)-V(X,Z)\dif \alpha(Y)\right)f_1,$$
and
$$\nabla_XV(Y,Z)-\nabla_YV(X,Z)+V(T(X,Y),Z)=0.$$
\EPf

\section{Curvature of transverse curves in the surface $S$}

Let be $\gamma\! :[a,b]\subset\mathbb R\rightarrow S$ a differentiable curve such that $\gamma'(t)$ is transversal, i.e., $f^2(\gamma'(t))\neq 0$ for all $t\in [a,b]$. Let be $T$ defined by
$$T(t)=\frac{1}{|f^2(\gamma'(t))|}\gamma'(t),$$
the unitary tangent field along $\gamma$. As $f^2(T(t))=\pm 1$, then
$$\nabla_Tf^2(T)+f^2(\nabla_TT)=0,$$
and as $\nabla f^2=0$, we know that 
$\nabla_TT$ is a multiple of $f_1$. We write
$$\nabla_TT=kN,$$
where the vector field $N=\epsilon f_1$ on $\gamma$, and $\epsilon=+1$ if $\{T,f_1\}$ is positively oriented and $\epsilon=-1$, otherwise.
Observe that  $\epsilon f^2(T)< 0$. The function $k\!:[a,b]\rightarrow \mathbb R$ is the \emph{curvature} of $\gamma$.
\begin{dfn}
The function $k=\langle\nabla_TT,N\rangle $ is the \emph{curvature} of the transverse curve $\gamma$.
\end{dfn}

\begin{prop}\label{fork}
The curvature $k$ is given by
$$k=\frac{\epsilon}{f^2(\gamma')}\left(\frac{d}{dt}\frac{f^1(\gamma')}{f^2(\gamma')}+A\,d\alpha(\gamma')\right).$$
\end{prop}
\Pf It follows from the definition that $k=\epsilon f^1(\nabla_TT)$, so $$k=\frac{\epsilon}{|f^2(\gamma')|}f^1(\nabla_{\gamma'}(\frac{1}{|f^2(\gamma')|}\gamma'))=\frac{\epsilon}{|f^2(\gamma')|}\left(\nabla_{\gamma'}(f^1(\frac{1}{|f^2(\gamma')|}\gamma'))-
(\nabla_{\gamma'}f^1)(\frac{1}{|f^2(\gamma')|}\gamma')\right),$$
and the proposition follows.\EPf

\section{Gauss-Bonnet theorem}

In this section, let be $R\subset S$ a fundamental set, and $c$ a fundamental $2$-chain such that $|c|=R$. The oriented curve $\gamma=\partial c$ is the bounding curve of $R$. The curve $\gamma$ is piecewise differentiable, and composed of differentiable curves $\gamma_j\!:[s_j,s_{j+1}]\rightarrow S$, $j=1,\ldots,r$, with $\gamma_j(s_{j+1})=\gamma_{j+1}(s_{j+1})$, for $j=1,\ldots,r-1$ and $\gamma_1(s_1)=\gamma_r(s_{r+1})$. We define the corner area at the \emph{vertices} $\gamma_j(s_{j+1})$ as
$ca_j=\ca(\gamma_{j}'(s_{j+1}),\gamma_{j+1}'(s_{j+1}))$, $j=1,\ldots,r-1$ and $ca_r=\ca(\gamma_{r}'(s_{r+1}),\gamma_{1}'(s_{1}))$.
\begin{thm}(Gauss-Bonnet formula) Let ${R}$ be contained in a coordinate domain $U$ of $S$, let the bounding curve $\gamma$ of $R$ be a simple closed transverse curve, and let $ca_1,\ldots,ca_r$ be the exterior corner areas of $\gamma$. Then
$$\int_{\gamma}k+\sum_{j=1}^rca_j+\int_RK=0,$$
where $k$ is the curvature function on $\gamma$ and $K$ is the  scalar curvature function on $R$.
\end{thm}
\Pf Let $\gamma_1,\ldots,\gamma_r$ be the $C^{\infty}$ pieces of $\gamma$ with $\gamma_j$ defined on the interval $[s_j,s_{j+1}]$, with $\gamma_j(s_{j+1})=\gamma_{j+1}(s_{j+1})$, for $j=1,\ldots,r-1$, and $\gamma_r(s_{r+1})=\gamma_1(s_1)$. Let be $\ca_j=\ca(\gamma_j'(s_j+1),\gamma_{j+1}'(s_j+1))$, for $j=1,\dots,r-1$ and $\ca_r=\ca(\gamma_r'(s_{r+1}),\gamma_1'(s_1))$. In each $C^\infty$ piece of $\gamma$ we have the positive orientation $T$ and the curvature $\nabla_TT=\epsilon kf_1$. Then from (\ref{K}), Propositions \ref{pnabla} and \ref{fork}, and
$$f_2=\frac{1}{f^2(\gamma_j')}\gamma_j'-\frac{f^1(\gamma_j')}{f^2(\gamma_j')}f_1=-\epsilon T-\frac{f^1(\gamma_j')}{f^2(\gamma_j')}f_1,$$
since $\epsilon=-\frac{|f^2(\gamma_j')|}{f^2(\gamma_j')}$, we obtain
$$\begin{aligned}
\int_RK&=\int_c Kf^1\wedge f^2=\int_c \dif A\wedge \dif \alpha(f_1,f_2)f^1\wedge f^2=\int_c \dif A\wedge \dif \alpha=\int_{\partial c}A\, \dif \alpha\\
&=\sum_{j=1}^r\int_{[s_j,s_{j+1}]}A\, \dif \alpha(\gamma_j')\\
&=\sum_{j=1}^r\int_{[s_j,s_{j+1}]}\left(\epsilon f^2(\gamma_j')k -\frac{\dif }{\dif t}\frac{f^1(\gamma_j')}{f^2(\gamma_j')}\right)\\
&=\int_{\partial c}- k|f^2| -\sum_{j=1}^r(\frac{f^1(\gamma_j'(s_{j+1}))}{f^2(\gamma_j'(s_{j+1}))}-\frac{f^1(\gamma_j'(s_{j}))}{f^2(\gamma_j'(s_{j}))})\vspace{.3cm}\\
&=-\int_{\gamma}k-\sum_{j=1}^r\frac{f^1(\gamma_j'(s_{j+1}))}{f^2(\gamma_j'(s_{j+1}))}+\sum_{j=0}^{r-1}\frac{f^1(\gamma_{j+1}'(s_{j+1}))}{f^2(\gamma_{j+1}'(s_{j+1}))}\vspace{.3cm}\\
&=-\int_{\gamma}k+\sum_{j=1}^{r-1}(\frac{f^1(\gamma_{j+1}'(s_{j+1}))}{f^2(\gamma_{j+1}'(s_{j+1}))}-\frac{f^1(\gamma_{j}'(s_{j+1}))}{f^2(\gamma_{j}'(s_{j+1}))})+
(\frac{f^1(\gamma_1'(s_{1}))}{f^2(\gamma_1'(s_{1}))}-\frac{f^1(\gamma_r'(s_{r+1}))}{f^2(\gamma_r'(s_{r+1}))})\vspace{.3cm}\\
&=-\int_{\gamma}k+\sum_{j=1}^{r-1}\ca(\gamma_{j+1}'(s_{j+1}),\gamma_{j}'(s_{j+1}))+\ca(\gamma_{1}'(s_{1}),\gamma_{r}'(s_{r+1}))\vspace{.3cm}\\
&=-\int_{\gamma}k-\sum_{j=1}^{r}\ca_j.
\end{aligned}
$$
If the surface $S$ is compact and oriented, then there exists a characteristic no-null vector field on $S$, therefore $S$ is diffeomorphic to a torus. In this case, we obtain the corollary:

\begin{cor}
Suppose $S$ is a differentiable compact surface in $\mathbb H^1$ with $\Sigma=\emptyset$. Then
$$\int_SK=0.$$
\end{cor}

\Pf In fact, we can triangulate $S$ by a finite number  of triangles $\Delta_i$, $i=1,\ldots,s,$ such that the boundary of each $\Delta_i$ is composed by transverse curves. As the triangles are positively oriented, then 
$$\int_SK=\sum_{i=1}^s\int_{\Delta_i}K=-\sum_{i=1}^s\int_{\partial\Delta_i}k-\sum_{i=1}^s\sum_{r=1}^3\mbox{ca}_{ir},$$
where $\partial\Delta_i$ is the boundary of $\Delta_i$ positively oriented, and $\mbox{ca}_{ir}$, $r\in\{1,2,3\}$, are the corner areas at each vertex of $\Delta_i$. If $\Delta_i$ and $\Delta_l$ have sides $\Delta_{iu}$ and $\Delta_{lv}$ in common, they have opposite orientations, so $\int_{\Delta_{iu}}k+\int_{\Delta_{lv}}k=0$; therefore, $\sum_{i=1}^s\int_{\partial\Delta_i}k=0$. In the same way, at a common vertex, the corner areas sum null, so $\sum_{i=1}^s\sum_{r=1}^3\mbox{ca}_{ir}=0$, and the proposition is proved.\EPf

\end{document}